\documentclass[a4paper,11pt]{article}

\newtheorem{theorem}{Theorem}
\newtheorem{proposition}{Proposition}
\newtheorem{lemma}{Lemma}
\newtheorem{corollary}{Corollary}
\newtheorem{remark}{Remark}

\newtheorem{example}{Example}
\newtheorem{assumption}{Assumption}

\usepackage{latexsym}
\usepackage{amsmath, amssymb}  
\usepackage{type1cm} 
\usepackage{bm}

\title{Global vs blow-up solutions and optimal threshold for hyperbolic ODEs with possibly singular nonlinearities}
\footnotetext{
(*) Corresponding author: \texttt{miyasita.t@yamato-u.ac.jp}}
\author{Daniele CASSANI \and Tosiya MIYASITA*}
\date{\today}
\begin{document}

\maketitle

\begin{abstract}
\noindent We consider a hyperbolic ordinary differential equation 
perturbed by a nonlinearity which can be singular at a point and in particular this includes MEMS type equations.
We first study qualitative properties of the solution to the stationary problem. Then, 
for small value of the perturbation parameter as well as initial value,
we establish the existence of a global solution by means of the Lyapunov function
and we show that the omega limit set consists of a solution to the stationary problem.
For strong perturbation or large initial value,
we show that the solution blows up.
Finally,
we discuss the relationship between upper bounds of the perturbation parameter 
for the existence of time-dependent and stationary solutions, for which we establish an optimal threshold. 
\end{abstract}

\section{Introduction}
\label{sec:intro}

In this paper, we study the following ordinary differential equation$:$
\begin{equation}
\left\{\begin{array}{ll}
u_{tt} + \alpha f(u_{t}) + \beta u = \lambda g(u)
& \mbox{for $t>0$}, \\
u(0) =  u_{0}, \\
u_{t}(0) = v_{0},
\end{array}\right.
\label{eqn:ODE-hyperbolic}
\end{equation}
where $\lambda>0$, $\alpha \geq 0$, $\beta > 0$, $u_{0} \in {\mathbb R}$ 
and 
$v_{0} \in {\mathbb R}$.
Under appropriate assumptions 
on $f$ and $g$,
we discuss the global existence and blow-up 
of the solutions 
to \eqref{eqn:ODE-hyperbolic}.
For $\alpha \geq 0$ 
and 
the initial value $\left( u_{0}, v_{0} \right) \in {\mathbb R}^2$,
we denote 
by $\lambda^* \left( u_{0}, v_{0} \right) (\alpha)$ 
the dynamical threshold for the existence 
of a global solution 
of \eqref{eqn:ODE-hyperbolic}.
Namely,
the solution of \eqref{eqn:ODE-hyperbolic} exists globally in time 
for $0 < \lambda < \lambda^* \left( u_{0}, v_{0} \right) (\alpha)$
and 
blows up 
for $\lambda > \lambda^* \left( u_{0}, v_{0} \right) (\alpha)$.
At the same time, let us denote by $\lambda^*$ the stationary threshold, see Theorem \ref{thm:stationary-sol} below, 
for the existence of solutions to 
\begin{equation}
\beta \phi = \lambda g ( \phi ).
\label{eqn:ODE-stationary}
\end{equation}
From the point of view of applications to Micro Electro Mechanical Systems, the value $\lambda^*$ plays an important role as it is connected to the so-called ``pull-in instability'', see \cite{CFT, CKL} and references therein. 
\newpage 

\begin{assumption}
We impose the following assumptions on the damping term $f$:
\begin{enumerate}
\item[$(f1)$] $f\in C^1({\mathbb R})$;
\item[$(f2)$] $f(0) = 0$;
\item[$(f3)$] $f(v)v>0$ for $v \in {\mathbb R} \setminus \left\{ 0 \right\}$;
\item[$(f4)$] there exist $\eta>0$ 
and 
$\theta \geq 0$ 
such that 
$\left\vert f(v) \right\vert \leq \eta \left\vert v \right\vert^{\theta+1}$
for $v \in {\mathbb R}$.
\end{enumerate}
\label{assumption:f}
\end{assumption}
\begin{remark} 
Note that under these assumptions one has $f'(0) \geq 0$.
\end{remark}

\begin{assumption}
Let $b \in (0, +\infty]$ and $I=(- \infty, b)$. We impose the following assumptions on the nonlinearity $g$:
\begin{enumerate}
\item[$(g1)$] $g(u) \in C^2(I)$;
\item[$(g2)$] $g(u)>0$ for $u \in I$.
Moreover, we normalize $g(u)$ by $g(0)=1$;
\item[$(g3)$] $g'(0) \geq 0$;
\item[$(g4)$] $g''(u) > 0$ for $u \in I$;
\item[$(g5)$] If $I= {\mathbb R}$,
for sufficiently large $R > 0$,
there exists $\tau>0$ such that $g''(u) > \tau$ for $u > R$.\\

If $I=(- \infty, b)$, $b < +\infty$,
the following hold:
\begin{equation*}
\lim_{u \to b^-}g(u)
= \lim_{u \to b^-}g'(u)
= \lim_{u \to b^-}g''(u)
= +\infty,
\end{equation*}
\begin{equation*}
\lim_{u \to b^-} \int_{0}^{u}g(s) \, ds = + \infty
\quad \text{and} \quad
\lim_{u \to b^-} \left( \int_{0}^{u}g(s) \, ds - \frac{u}{2}g(u) \right) <0.
\end{equation*}
\end{enumerate}
\label{assumption:g}
\end{assumption}
\begin{remark} 
Note that under these assumptions solutions to \eqref{eqn:ODE-stationary} are positive.
\end{remark}
So far, a huge literature has been devoted to \eqref{eqn:ODE-hyperbolic} and related problems from the theoretical point of view as well as from the point of view of applications, see for instance \cite{VBM, KLNT, GK} and references therein. 
 
\noindent In \cite{MR3592652},
Flores studies the following problem 
\begin{equation*}
\left\{\begin{array}{ll}
u_{tt} + \alpha u_{t} + u = \frac{\lambda}{\left( 1-u \right)^2},
& \mbox{for $t>0$}, \\
u(0) =  u_{0} \in [0,1), \\
u_{t}(0) = v_{0}.
\end{array}\right.
\label{eqn:ODE-flores}
\end{equation*}
He proves that $0 < \lambda \left( 0,0 \right) (\alpha) < \lambda^*$ 
for $\alpha \geq 0$
and that 
$\lim_{\alpha \to +\infty} \lambda \left( 0,0 \right) (\alpha) = \lambda^*$.
In \cite{MR2763360},
Haraux considers the following 
\begin{equation*}
u_{tt} 
+ c \left\vert u_{t} \right\vert^{\alpha} u_{t} 
+ \left\vert u \right\vert^{\beta} u
= 0,
\qquad
\text{for $c>0$, $\alpha>0$, $\beta>0$}
\label{eqn:MR2763360-homo}
\end{equation*}
and
studies the existence of sign-changing solutions,
number of the zeros 
and
the decay estimates for solutions
depending on the value
of the parameters  $c$, $\alpha$, $\beta$.
Moreover,
the results are generalized to 
\begin{equation*}
u_{tt} 
+ c \left\vert u_{t} \right\vert^{\alpha} u_{t} 
+ \left\vert u \right\vert^{\beta} u
= f(t),
\qquad
\text{for $c>0$, $\alpha>0$, $\beta>0$}
\label{eqn:MR2763360-inhomo}
\end{equation*}
for a continuous function $f(t)$
with some decaying properties.
In \cite{MR1608005},
Souplet studies the following backward equation 
\begin{equation*}
u_{tt} 
- \left\vert u_{t} \right\vert^{\alpha} u_{t}
+ \left\vert u \right\vert^{\beta} u
= 0
\qquad
\text{for $\alpha>0$, $\beta>0$}
\label{eqn:MR1608005-homo}
\end{equation*}
and
proves  the existence of unbounded global solutions,
unbounded oscillatory solutions, as well as their blow-up rate
and
asymptotic behaviour.
In \cite{MR1000727, MR1270119},
the authors consider
\begin{equation*}
u_{tt} 
- \left\vert u_{t} \right\vert^{\alpha}
+ \lambda \left\vert u \right\vert^{\beta}
= 0,
\qquad
\text{for $\lambda>0$, $\alpha>0$, $\beta>0$}.
\label{eqn:MR1270119-homo}
\end{equation*}
They investigate for which value 
of parameters one has existence of solutions,
derive their asymptotic behaviour
and
classify the ground state solutions
according to the value
of parameters.\\
Here we aim at extending the result in \cite{MR3592652} 
to more general nonlinearities $f$ and $g$.
Henceforth,
we consider \eqref{eqn:ODE-hyperbolic} 
and 
\eqref{eqn:ODE-stationary}
under Assumptions \ref{assumption:f} 
and 
\ref{assumption:g}
unless otherwise stated.
We consider first the stationary problem \eqref{eqn:ODE-stationary} for which 
we obtain the bifurcation diagram of the solution set $\{ \left( \lambda, \phi \right) \}$.

\begin{theorem}
There exists a unique $p>0$ such that 
\begin{equation*}
g(p) - p g'(p)
= 0.
\label{eqn:p}
\end{equation*}
Let 
\begin{equation*}
\lambda^*
:= \frac{\beta p}{g(p)}.
\label{eqn:static-lambda}
\end{equation*}
The following hold:
\begin{enumerate}
\item[i)] For any $\lambda < \lambda^*$,
there are two solutions $\phi_{1} = \phi_{1} ( \lambda )$ 
and $\phi_{2} = \phi_{2} ( \lambda )$ 
of \eqref{eqn:ODE-stationary} with
$0 < \phi_{1} < p < \phi_{2}<b$
and
\begin{equation*}
\lim_{\lambda \searrow 0} \left( \phi_{1}, \phi_{2} \right)
= \left( 0, b \right) 
\quad \text{and} \quad
\lim_{\lambda \nearrow \lambda^*} \left( \phi_{1}, \phi_{2} \right)
=: \left( \phi_{1}^*, \phi_{2}^* \right)
= \left( p, p \right) \ .
\end{equation*}

Moreover,
$\phi_{1}( \lambda )$ and $\phi_{2}( \lambda )$ are respectively increasing and decreasing
with respect to  $\lambda$; 
\item[ii)] For $\lambda = \lambda^*$,
there exists a unique solution $\phi_{1} = \phi_{2} = p$ of \eqref{eqn:ODE-stationary};
\item[iii)] For $\lambda > \lambda^*$,
do not exist solutions to \eqref{eqn:ODE-stationary}.
\end{enumerate}
\label{thm:stationary-sol}
\end{theorem}

\noindent In order to investigate the dynamical behaviour,
we find the solution $( {\overline \lambda}, \overline{\phi_{2}} )$ 
with the following properties:

\begin{theorem}
There exists a unique solution  
$( \lambda, \phi_{2}(\lambda) ) = ( {\overline \lambda}, \overline{\phi_{2}} )$ 
of \eqref{eqn:ODE-stationary}
satisfying the following:
\begin{equation*}
0 < \overline{\lambda} < \lambda^*,
\quad
p < \overline{\phi_{2}}<b
\quad \text{and} \quad
{\overline \lambda}
= \frac{\beta}{2} \frac{\overline{\phi_{2}}^2}{\int_{0}^{\overline{\phi_{2}}} g(s) \, ds}.
\label{eqn:dyn-lambda}
\end{equation*}
\label{thm:stationary-dyn-sol}
\end{theorem}

\noindent Next,
on the one hand we consider the time-dependent equation \eqref{eqn:ODE-hyperbolic}
and 
derive the conditions 
for the existence of global bounded solutions.

\begin{theorem}
For $\lambda < \lambda^*$ 
and 
$\alpha \geq 0$,
let 
\begin{equation*}
D_{0}:=\left\{ \left( u, v \right) \in {\mathbb R}^2 \mid
l(\lambda) < u < \phi_{2}(\lambda),
\
v^2 < J_{\lambda}(\phi_{2}) - J_{\lambda}(u_{0}) \right\},
\label{eqn:initial-sp}
\end{equation*}
where 
$l(\lambda)$,
$J_{\lambda}(\phi_{2})$
and
$J_{\lambda}(u_{0})$
are constants 
depending only on $\lambda$, $\beta$, $g$ and $u_{0}$.
If $\left( u_{0}, v_{0} \right) \in D_{0}$,
then \eqref{eqn:ODE-hyperbolic} has a unique global solution 
$u \in W^{2,\infty} \left( [0,\infty) \right)$.
Moreover,
if $\alpha > 0$,
then the following holds
\begin{equation}
\lim_{t \to +\infty} 
\Big( \left\vert u(t) - \phi_{1}(\lambda) \right\vert 
+ \left\vert u_{t}(t) \right\vert \Big) 
= 0\ .
\label{eqn:convergence}
\end{equation}
\label{thm:ex-global-sol}
\end{theorem}

\noindent On the other hand,
we obtain the following blow-up result for $\lambda > \lambda^*$ or $\phi_{2}<u_{0}<b$.
Let $T_{\infty} \in (0, +\infty]$, be the maximal time of existence of a solution to \eqref{eqn:ODE-hyperbolic}. Then, we have the following 

\begin{theorem}
Let $\lambda > \lambda^*$,
$\alpha \geq 0$ and $u_{0},v_{0}  \geq 0$.
The solution $(u, u_{t})$ of \eqref{eqn:ODE-hyperbolic} 
blows up to $( +\infty, +\infty )$ for $b=+\infty$
and 
quenches to $( b, +\infty )$ for $0 < b <+\infty$
as $t \to T_{\infty}$,
where
$T_{\infty} < +\infty$
if $\alpha = 0$
or 
$0 < b <+\infty$.
\label{thm:blowup-large-para}
\end{theorem}

\begin{theorem}
Let $0 < \lambda \leq \lambda^*$ 
and 
$\alpha \geq 0$.
For any $\phi_{2}(\lambda) < u_{0} < b$ and $v_{0} \geq 0$,
the same conclusion of Theorem \ref{thm:blowup-large-para} holds.
\label{thm:blowup-large-IV}

\end{theorem}

\noindent The main result of this paper is concerned with establishing an optimal dynamical threshold, indeed we have 

\begin{theorem}
Let $f(v)=v$. The function 
$\lambda \left( 0,0 \right) (\alpha)$ is continuous and monotone increasing 
with respect to $\alpha \geq 0$ and satisfies:
\begin{itemize}
\item[$i)$] $\overline{\lambda} \leq \lambda \left( 0,0 \right) (\alpha) < \lambda^*$ 
for $\alpha \geq 0$;
\item[$ii)$] $\lim_{\alpha \to 0^+} \lambda \left( 0,0 \right) (\alpha) = \overline{\lambda}$;
\item[$iii)$] $\lim_{\alpha \to +\infty} \lambda \left( 0,0 \right) (\alpha) = \lambda^*$.
\end{itemize}
\label{thm:dynamical-threshold}
\end{theorem}

\noindent This paper is organized as follows$:$
In Section \ref{sec:stationary problem},
we consider the stationary problem.
Thanks to the assumptions on $g$,
we get at most two positive solutions.
In Section \ref{sec:dynamical system}, we
settle preliminary lemmas, involving energy and dynamical system, 
in order to investigate the dynamical behaviour.
On the one hand, in Section \ref{sec:time-depending problem},
we establish the existence of a global solution and periodic orbit under some appropriate conditions.
On the other hand,
for large values of the perturbation parameter or large initial values,
we prove that the solution blows up.
In Section \ref{sec:properties of the orbit},
we discuss qualitative properties of the orbit such as openness, monotonicity and continuity.
In Section \ref{sec:proof-of-main-TH},
we prove our main result, namely Theorem \ref{thm:dynamical-threshold}.

\section{The stationary problem}
\label{sec:stationary problem}

Here we study the solution set of the function equation \eqref{eqn:ODE-stationary}
and obtain the upper bound $\lambda^*$ for the existence of solutions.
We regard $\phi = \phi ( \lambda )$ as a function of the parameter $\lambda$.

\medskip

\noindent {\it Proof of Theorem \ref{thm:stationary-sol}.}
Since solutions are positive,
we consider  in \eqref{eqn:ODE-stationary} $\phi \geq 0$.
By Assumption \ref{assumption:g},
\begin{equation*}
F(u):=\frac{u}{g(u)}
\label{eqn:F}
\end{equation*}
is well-defined for $u \geq 0$
and the following holds 
\begin{equation*}
F'(u)
= \frac{1}{g(u)^2}
\left( g(u) - u g'(u) \right)
\label{eqn:F'}\ .
\end{equation*}
From  $F(0)=F(b)=0$ and $F'(0)>0$,
there exists $p>0$ such that $F'(p)=0$.
Next we show the uniqueness of such $p$.
Set $G(u) := g(u) - u g'(u)$,
and consider the sign of $G(u)$.
We have $G(0)=1$,
\begin{equation*}
G'(u)
= - u g''(u)
< 0
\label{eqn:G'}
\end{equation*}
for all $0<u<b$,
$G'(0)=0$
and $G'(b) = -\infty$,
which implies that $G(u)>0$ for $u<p$ and $G(u)<0$ 
for $u>p$
and thus
such $p$ is uniquely determined.
The statement follows by drawing the graph 
of $y=\beta F(u)$ 
and 
$y=\lambda$.
\hfill $\Box$

\vspace{1pc}

\noindent Define
\begin{equation*}
H_{\lambda}(u):=\beta u - \lambda g(u).
\label{eqn:H}
\end{equation*}
Then
$\beta F(u) =\lambda$
is equivalent to
$H_{\lambda}(u)=0$ and the sign of $H_{\lambda}(u)$ 
is given by the value of $\lambda$.

\begin{corollary}
For $\lambda < \lambda^*$,
we have:
\begin{itemize}
\item[i)] $H_{\lambda}(u) < 0$, 
for $u < \phi_{1}$, $\phi_{2}<u<b$;
\item[ii)] $H_{\lambda}(u) = 0$, 
for $u = \phi_{1}$, $u = \phi_{2}$;
\item[iii)] $H_{\lambda}(u) > 0$, 
for $\phi_{1} < u < \phi_{2}$.
\end{itemize}
Furthermore,
$H_{\lambda}(u)$ is increasing in $(- \infty, q)$ and decreasing 
in $(q, b)$, where $q$ satisfies $\beta = \lambda g'(q)$ 
and 
$p < q < \phi_{2}<b$.
\label{coro:H-sign}
\end{corollary}
{\it Proof.}
Noting that $H_{\lambda}(\phi_{1}) = H_{\lambda}(\phi_{2}) = 0$, $H_{\lambda}'(u)
= \beta - \lambda g'(u)$ and that $H_{\lambda}''(u)= - \lambda g''(u)< 0$, 
we find $q \in (\phi_{1}, \phi_{2})$ uniquely satisfying $H_{\lambda}'(q)=0$.
Moreover, from 
\begin{equation*}
H_{\lambda}'(p)
= \beta - \lambda g'(p)
= \lambda^* \frac{g(p)}{p} - \lambda \frac{g(p)}{p}
= \left( \lambda^* - \lambda \right) \frac{g(p)}{p}
> 0,
\label{eqn:H'-p}
\end{equation*}
one has $p<q$.
\hfill $\Box$

\begin{remark}
For $\lambda = \lambda^*$,
we have $\phi_{1}(\lambda^*) = \phi_{2}(\lambda^*) = p$ and $H_{\lambda^*}(u) < 0$, 
for $u \neq p$, $H_{\lambda^*}(p) = 0$. Furthermore,
$H_{\lambda^*}(u)$ is increasing in $(- \infty, p)$ and decreasing in $(p,b)$.
\label{coro:H^*-sign}
\end{remark}

\begin{corollary}
For $\lambda > \lambda^*$,
we have 
\begin{equation*}
H_{\lambda}(u) < 0
\quad \text{for $u<b$}.
\end{equation*}
In particular,
there exists $\xi > 0$ such that 
\begin{equation*}
H_{\lambda}(u) \leq - \xi < 0
\end{equation*}
for all $u \geq 0$,
where $\xi$ depends only on $\lambda$, $\beta$ and $g$.
\label{coro:H-large-sign}
\end{corollary}
{\it Proof.}
Set $\xi = \lambda - \lambda^*$,
to have
\begin{equation*}
\xi
= \lambda - \beta F(p)
\leq \lambda - \beta F(u)
= - \frac{1}{g(u)} H_{\lambda}(u)
\leq - H_{\lambda}(u)
\end{equation*}
as
$1 = g(0) \leq g(u)$ holds for all $u \geq 0$.
\hfill $\Box$

\vspace{1pc}

\begin{lemma}
For $\lambda < \lambda^*$, one has 
\begin{equation*}
H_{\lambda}'(\phi_{1}(\lambda))
= \sup_{\phi_{1} \leq u \leq \phi_{2}} \frac{H_{\lambda}(u)}{u - \phi_{1}(\lambda)} \ .
\label{eqn:condition-heteroclinic}
\end{equation*}
\label{coro:heteroclinic-criterion}
\end{lemma}
{\it Proof.}
Noting that $H_{\lambda}(\phi_{1})=0$,
we apply the mean value theorem to obtain the conclusion, as $H_{\lambda}''(u)<0$ for $\phi_{1} \leq u \leq \phi_{2}$.
\hfill $\Box$

\medskip

\noindent {\it Proof of Theorem \ref{thm:stationary-dyn-sol}.}
Let us define
\begin{equation*}
I(u):= 2 \int_{0}^{u} g(s) \, ds - u g(u), \quad 0\leq u<b\ .
\label{eqn:I}
\end{equation*}
We \textit{claim} that $I(u_{0}) = 0$ for some $p<u_{0} <b$. As a consequence, we find a solution $( \lambda, \phi_{2} ) = ( \beta u_{0} /g( u_{0} ), u_{0} )$ 
of \eqref{eqn:ODE-stationary}
by Theorem \ref{thm:stationary-sol}
and
\begin{equation*}
\lambda
= \frac{\beta u_{0}}{g( u_{0} )}
= \frac{\beta u_{0}^2}{2 \int_{0}^{u_{0}} g(s) \, ds}
\in (0, \lambda^*).
\label{eqn:I-zero}
\end{equation*}
Let us prove the claim. First we have
\begin{eqnarray*}
I( p )
&=& 2 \int_{0}^{p} g(s) \, ds - p g(p) \\
&=& 2 \int_{0}^{p} g(s) \, ds - \frac{\beta}{\lambda^*} p^2 \\
&=& \frac{2}{\lambda^*} \int_{0}^{p} \left( \lambda^* g(s) - \beta s \right) \, ds \\
&=& - \frac{2}{\lambda^*} \int_{0}^{p} H_{\lambda^*}(s) \, ds>0 \\
\label{eqn:I^*-est}
\end{eqnarray*}
by Remark \ref{coro:H^*-sign}.
Observe that $I'(u) = G(u) < 0$ 
and that 
$I''(u) = G'(u) < 0$ 
for $p<u<b$,
where $G$ is defined as in the proof of Theorem \ref{thm:stationary-sol}.
By the monotonicity of $I$,
we need only to show that $I(b)<0$.
In the case of $b < +\infty$,
$I(b)<0$ follows from Assumption \ref{assumption:g}.
In the case $b=+\infty$, 
since
$I''(+\infty) = G'(+\infty) = -\infty$ holds,
we have $I'(+\infty) = G(+\infty) = -\infty$
and finally 
$I(+\infty) = -\infty$,
which proves the theorem.
\hfill $\Box$

\medskip

\noindent Let us give a few examples where $( \lambda^*, p )$ 
and 
$( {\overline \lambda}, \overline {\phi_{2}} )$ are explicitly known. 

\begin{example}
If $\beta=1$ and $g(u) = \left( 1+u \right)^{2k}$ for $k \in {\mathbb N}$,
Assumption \ref{assumption:g} is satisfied with $b=+\infty$.
In particular,
in the case of $k=1$,
we have
\begin{equation*}
\left( \lambda^*, p \right) = \left( \frac{1}{4}, 1 \right)
\quad \text{and} \quad
( {\overline \lambda}, \overline {\phi_{2}} )
= \left( \sqrt{3} - \frac{3}{2}, \sqrt{3} \right).
\label{eqn:dyn-power}
\end{equation*}
\label{theom:ex-power}
\end{example}

\begin{example}
If $\beta=1$ and $g(u) = e^u$,
Assumption \ref{assumption:g} is satisfied with $b=+\infty$.
We have
\begin{equation*}
\left( \lambda^*, p \right) = \left( \frac{1}{e}, 1 \right)
\label{eqn:ext-exp}
\end{equation*}
and $\left( {\overline \lambda}, \overline {\phi_{2}} \right)$ satisfies
\begin{equation*}
\overline {\phi_{2}} e^{\overline {\phi_{2}}} - 2 e^{\overline {\phi_{2}}} + 2
= 0,
\quad
{\overline \lambda}
= \overline {\phi_{2}} e^{-\overline {\phi_{2}}},
\quad 
\overline{\phi_{2}} \in (1, 2).
\label{eqn:dyn-exp}
\end{equation*}
\label{theom:ex-exp}
\end{example}

\begin{example}
If $\beta=1$ and $g(u) = 1+u^{2k}$ for $k \in {\mathbb N}$,
Assumption \ref{assumption:g} is satisfied with $b=+\infty$.
In particular,
in the case $k=1$,
we have
\begin{equation*}
\left( \lambda^*, p \right) = \left( \frac{1}{2}, 1 \right)
\quad \text{and} \quad
( {\overline \lambda}, \overline {\phi_{2}} )
= \left( \frac{\sqrt{3}}{4}, \sqrt{3} \right). 
\label{eqn:dyn-poly}
\end{equation*}
\label{theom:ex-poly}
\end{example}

\begin{example}
If $\beta=1$ and $g(u) = 1 / \left( 1-u \right)^p$ for $p>1$,
$g(u)$ has a singularity
at $u=1$.
However,
Assumption \ref{assumption:g} is satisfied
with $b=1$.
In particular,
in the case of $p=2$ we have
\begin{equation*}
\left( \lambda^*, p \right) = \left( \frac{4}{27}, \frac{1}{3} \right)
\quad \text{and} \quad
( {\overline \lambda}, \overline {\phi_{2}} )
= \left( \frac{1}{8}, \frac{1}{2} \right). 
\label{eqn:dyn-flores}
\end{equation*}
\label{theom:ex-flores}

\end{example}

\section{A dynamical system}
\label{sec:dynamical system}

Consider the following energy functional
\begin{equation*}
E_{\lambda}( u, v ):=\frac{\beta}{2} u^2 + \frac{1}{2} v^2 - \lambda \int_{0}^{u} g(s) \, ds.
\end{equation*}
Then,
$E_{\lambda}( u(t), u_{t}(t) )$ turns out to be the Lyapunov function 
for \eqref{eqn:ODE-hyperbolic}.
In fact, we have
\begin{equation*}
\frac{d}{dt} E_{\lambda}( u(t), u_{t}(t) )
= -\alpha f(u_{t}(t)) u_{t}(t)
\leq 0,
\end{equation*}
which yields
\begin{equation}
E_{\lambda}( u(t), u_{t}(t) ) + \alpha \int_{0}^{t} f(u_{t}(r))u_{t}(r) \, dr
= E_{\lambda}( u_{0}, v_{0} ) .
\label{eqn:integration-lyaunov}
\end{equation}
Hence,
by defining
\begin{equation*}
J_{\lambda}(u)
:= \frac{\beta}{2} u^2 - \lambda \int_{0}^{u} g(s) \, ds,
\label{eqn:J}
\end{equation*}
the energy inequality
\begin{equation}
J_{\lambda}(u)
\leq J_{\lambda}(u_{0}) + \frac{1}{2}v_{0}^2
\label{eqn:energy-ineq}
\end{equation}
holds by \eqref{eqn:integration-lyaunov}.
Every local solution satisfies \eqref{eqn:energy-ineq} as long as it exists.
To extend the solution globally in time,
we consider some properties of $J_{\lambda}(u)$ 
in the following two lemmas:

\begin{lemma}
For $\lambda < \lambda^*$,
$J_{\lambda}(u)$ has a local minimum
at $u=\phi_{1}(\lambda)$ and a local maximum $u=\phi_{2}(\lambda)$.
Moreover,
$J_{\lambda}(\phi_{1}(\lambda))<0$ holds.
\label{lem:local-min-max}
\end{lemma}
{\it Proof.}
Since $J_{\lambda}'(u) = \beta u - \lambda g(u) = H_{\lambda}(u)$ 
and 
$J_{\lambda}''(u) = H_{\lambda}'(u)$,
we have 
$J_{\lambda}'(\phi_{1}) = J_{\lambda}'(\phi_{2}) = 0$ 
and 
$J_{\lambda}''(\phi_{2}) < 0 < J_{\lambda}''(\phi_{1})$
by Corollary \ref{coro:H-sign}.
Then
$J_{\lambda}(\phi_{1}(\lambda))
< J_{\lambda} ( 0 )
= 0$
for all $\lambda \in (0, \lambda^*)$,
which completes the proof.
\hfill $\Box$

\vspace{1pc}

\begin{lemma}
The following hold:
\begin{enumerate}
\item[i)] $J_{\lambda}(\phi_{2}(\lambda)) >0 $ 
for $0 < \lambda < {\overline \lambda}$;
\item[ii)] $J_{\lambda}(\phi_{2}(\lambda)) =0 $ 
for $\lambda = {\overline \lambda}$;
\item[iii)] $J_{\lambda}(\phi_{2}(\lambda)) <0 $ 
for ${\overline \lambda} < \lambda < \lambda^*$.
\end{enumerate}
\label{lem:local-max}
\end{lemma}
{\it Proof.}
From
\begin{equation*}
\frac{d}{d \lambda} J_{\lambda}(\phi_{2}(\lambda))
= H_{\lambda} \left( \phi_{2}(\lambda) \right) \frac{d}{d \lambda} \phi_{2}(\lambda) 
- \int_{0}^{\phi_{2}(\lambda)} g(s) \, ds
= - \int_{0}^{\phi_{2}(\lambda)} g(s) \, ds
< 0
\label{eqn:local-min}
\end{equation*}
by Corollary \ref{coro:H-sign},
$J_{\lambda}(\phi_{2}(\lambda))$ is monotone decreasing in $\lambda$.
It follows from Theorem \ref{thm:stationary-dyn-sol} that 
$J_{\overline{\lambda}}(\phi_{2}(\overline{\lambda})) = 0$,
which yields the conclusions.
\hfill $\Box$

\vspace{1pc}

\noindent In the next section,
we consider dynamical properties
of the solution
of \eqref{eqn:ODE-hyperbolic},
which can be written in the following form 
\begin{equation}
\frac{d}{dt}
\left[ \begin{array}{cc}
u\\
u_{t}
\end{array} \right]
= \left[ \begin{array}{cc}
u_{t}\\
- \alpha f(u_{t}) - H_{\lambda}(u)
\end{array} \right].
\label{eqn:dyn-sys}
\end{equation}
Now
under Assumptions \ref{assumption:f} 
and 
\ref{assumption:g},
we obtain a local solution.
Next,
we establish the existence of a global solution 
exploiting Lemmas \ref{lem:local-min-max} and \ref{lem:local-max}.
For this purpose,
we consider the stability of the equilibrium point. 
At the equilibrium point
\begin{equation*}
\left( u, u_{t} \right)
= \left( \phi_{i}(\lambda), 0 \right)
\label{eqn:equilibrium-pt}
\end{equation*}
for $i=1$, $2$,
the linearized equation is given by
\begin{equation*}
\frac{d}{dt}
\left[ \begin{array}{cc}
U\\
V
\end{array} \right]
= \left[ \begin{array}{cc}
0&
1\\
-H_{\lambda}'(\phi_{i})&
-\alpha f'(0)
\end{array} \right]
\left[ \begin{array}{cc}
U\\
V
\end{array} \right]
\label{eqn:lineraized-eqn}
\end{equation*}
and
the eigenvalues  $\mu_{i}^{\pm}$
of coefficient matrix 
are given respectively by
\begin{equation*}
\mu_{i}^{\pm}
= \frac{- \alpha f'(0)  \pm \sqrt{\left\{ \alpha f'(0) \right\}^2 - 4 H_{\lambda}'(\phi_{i})}}{2}
\label{eqn:liearized-op}\ .
\end{equation*}
Along with Corollary \ref{coro:H-sign},
we have the following two lemmas$:$

\begin{lemma}
For $\lambda < \lambda^*$,
the equilibrium point
$\left( u, u_{t} \right) = \left( \phi_{1}(\lambda), 0 \right)$
is a stable focus 
for $0 < \alpha f'(0) < 2 \sqrt{H_{\lambda}'(\phi_{1} (\lambda) )}$,
a stable node 
for $\alpha f'(0) \geq 2 \sqrt{H_{\lambda}'(\phi_{1} (\lambda) )}$ and it is a centre
for $\alpha f'(0)=0$.
\label{lem:eqilibria-sink-centre}
\end{lemma}


\begin{lemma}
For $\lambda < \lambda^*$,
the equilibrium point
$\left( u, u_{t} \right) = \left( \phi_{2}(\lambda), 0 \right)$
is a saddle for all $\alpha \geq 0$.
\label{lem:eqilibria-saddle}
\end{lemma}

\section{The time-dependent problem}
\label{sec:time-depending problem}

In this Section,
for small parameters and small initial values,
we establish the existence of a global solution by means of the Lyapunov function method. 
For the dissipative case, that is $\alpha >0$,
we show that the global solution converges to the stationary solution.
For the conservative case, namely $\alpha = 0$,
we consider the periodic orbit starting at $(0,0)$.
Finally,
we show that the solution becomes unbounded
for large parameters as well as for large initial values.

\vspace{1pc}

\noindent {\it Proof of Theorem \ref{thm:ex-global-sol}.}
We define by $l(\lambda) \in ( -\infty, \phi_{1}(\lambda) )$ 
the point satisfying $J_{\lambda}(u) = J_{\lambda}( \phi_{2}(\lambda) )$.
By Lemma \ref{lem:local-max},
we have $l(\lambda)<0$, 
$l(\lambda)=0$ 
and 
$l(\lambda)>0$
for $\lambda \in (0, {\overline \lambda})$, 
$\lambda = {\overline \lambda}$ 
and 
$\lambda \in ( {\overline \lambda}, \lambda^*)$,
respectively. 
Note that
$J_{\lambda}(u_{0}) < J_{\lambda}(\phi_{2})$
holds
for
$l(\lambda) < u_{0} < \phi_{2}(\lambda)$
by Lemma \ref{lem:local-min-max}.
Then,
for $\left( u_{0}, v_{0} \right) \in D_{0}$,
\eqref{eqn:energy-ineq} yields
\begin{equation*}
J_{\lambda}(u)
< \frac{1}{2} \left( J_{\lambda}(\phi_{2}) +  J_{\lambda}(u_{0}) \right)
< J_{\lambda}(\phi_{2})
\label{eqn:energy-est}
\end{equation*}
and moreover
\begin{equation}
l(\lambda)
< u(t)
< J_{\lambda}^{-1} \left( \frac{J_{\lambda}(\phi_{2}) +  J_{\lambda}(u_{0})}{2} \right)
< \phi_{2},
\label{eqn:sol-est-1}
\end{equation}
where $J_{\lambda}^{-1}$ is the inverse function 
of $J_{\lambda}$ defined 
at $(\phi_{1}, \phi_{2})$.
Finally
\eqref{eqn:integration-lyaunov} implies that
\begin{equation*}
u_{t}^2
\leq 2 E_{\lambda}( u_{0}, v_{0} ) + 2 \lambda \int_{0}^{u} g(s) \, ds 
< 2 E_{\lambda}( u_{0}, v_{0} ) + 2 \lambda g(\phi_{2}) \phi_{2}
\label{eqn:sol-est-2}
\end{equation*}
and that $\left( u(t), u_{t}(t) \right)$ is uniformly bounded 
in ${\mathbb R}^2$ for all $t \geq 0$.
Therefore, 
since $u_{tt}$ is also bounded for $t \geq 0$ by \eqref{eqn:ODE-hyperbolic},
we have $u(t) \in W^{2, \infty}\left( [0,\infty) \right)$.
If $\alpha >0$,
the Lyapunov function $E_{\lambda} \left( u(t), u_{t}(t) \right)$ is strictly decreasing 
in $t$.
Thus
it follows 
from Theorem 5.1.8 and Corollary 8.5.1 
in \cite{MR3380967} 
that
the omega limit set $\omega \left( u_{0}, v_{0} \right)$ is connected in ${\mathbb R}^2$
and 
included 
in the solution set of \eqref{eqn:ODE-stationary}.
Since $\overline{\cup_{t \geq 0} \left( u(t), u_{t}(t) \right)} 
\not \in 
\left\{ \left( \phi_{2}, 0 \right) \right\}$
by \eqref{eqn:sol-est-1},
we obtain \eqref{eqn:convergence}.
\hfill $\Box$

\begin{remark}
Since $l (\lambda) < \phi_{1}$ 
and 
$J_{\lambda}(u) < J_{\lambda}(\phi_{2})$ 
for all $\lambda < \lambda^{*}$
and
$\phi_{1} < u < \phi_{2}$,
we have $(u_{0}, 0) \in D_{0}$ with $\phi_{1} < u_{0} < \phi_{2}$.
In other words,
the solution exists globally in time
for the initial value $(u_{0}, v_{0})$ with $\phi_{1} < u_{0} < \phi_{2}$ and $v_{0} = 0$.
\label{rem:ex-global}
\end{remark}

\vspace{1pc}

\noindent {\it Proof of Theorem \ref{thm:blowup-large-para}.} Let us divide the proof into four cases: 

\medskip

\noindent {\bf I. The case $\alpha = 0$ and $I=(-\infty, b)$.}\\
Integrating \eqref{eqn:ODE-hyperbolic}
and thanks to Corollary \ref{coro:H-large-sign},
we have
\begin{equation}
u_{t}
\geq v_{0} + \xi t
\geq 0
\quad 
\text{and} 
\quad
u
\geq u_{0} + v_{0}t + \frac{\xi}{2} t^2
\geq 0,
\label{eqn:blowup-02}
\end{equation}
which implies
that
$u(t)$ reaches  $b$
for finite $T_{\infty} < \infty$.
Then
\begin{equation*}
u_{t}^2 ( T_{\infty} )
= 2 E_{\lambda}( u_{0}, v_{0} ) + 2 \lambda \int_{0}^{b} g(s) \, ds 
- \beta b^2
= +\infty
\label{eqn:sol-est-2}
\end{equation*}
by \eqref{eqn:integration-lyaunov} and
Assumption \ref{assumption:g}.

\medskip

\noindent {\bf II. The case $\alpha = 0$ and $I={\mathbb R}$.}\\
Assume by contradiction
that
$T_{\infty} = \infty$.
By \eqref{eqn:blowup-02},
there exists sufficiently large $T>0$
such that
\begin{equation}\label{dceq1}
g''(u(t)) > \tau
\end{equation}
for all $t > T$, where $\tau$ is given 
in Assumption \ref{assumption:g}
and
$T$ can be taken 
as $T=\sqrt{2 \xi^{-1} R}$.
Then, 
integrating inequality \eqref{dceq1}
twice, with respect to $u$,
over $[u(T), u]$ we get
\begin{equation*}
g(u(t)) 
\geq g(u(T)) + \frac{\tau}{2} \left( u(t) - u(T) \right)^2
\end{equation*}
as $g'(u) \geq 0$
for $u \geq 0$.
Thus
we have
\begin{eqnarray}
u_{tt}(t)
&=& \lambda g(u(t)) - \beta u(t) \nonumber \\
&\geq& \lambda \left\{ g(u(T)) + \frac{\tau}{2} \left( u(t) - u(T) \right)^2 \right\}
- \beta u(t) \nonumber \\
&=& -H_{\lambda}(u(T)) 
+ \frac{\tau \lambda}{2} \left( u(t) - u(T) \right)^2 
- \beta \left( u(t) - u(T) \right) \nonumber \\
&>& \frac{\tau \lambda}{2} \left( u(t) - u(T) \right)^2 
- \beta \left( u(t) - u(T) \right)
\label{eqn:blowup-10}
\end{eqnarray}
for all $t > T$.
Since $u'(t) > 0$ holds
for all $t > T$
by \eqref{eqn:blowup-02},
we have
\begin{eqnarray*}
u_{t} u_{tt} 
&>& \frac{\tau  \lambda}{2} \left( u(t) - u(T) \right)^2 u_{t}
- \beta \left( u(t) - u(T) \right) u_{t}
\end{eqnarray*}
and
\begin{eqnarray*}
\left\{ u_{t}(t) \right\}^2
&>& \left\{ u_{t}(T) \right\}^2
+ \frac{\tau \lambda}{3} \left( u(t) - u(T) \right)^3 
- \beta \left( u(t) - u(T) \right)^2 \\
&>& \left\{ \frac{\tau \lambda}{3} \left( u(t) - u(T) \right) -  \beta \right\}
\left( u(t) - u(T) \right)^2.
\end{eqnarray*}
Take $T_{1} \in ( T, +\infty)$
such that
\begin{equation*}
\frac{\tau \lambda}{3} \left( u(t) - u(T) \right) -  \beta 
> \frac{\tau \lambda}{6} \left( u(t) - u(T) \right)
\end{equation*}
holds
for all $t > T_{1}$.
For instance,
we can take $T_{1}$ 
as follows  
\begin{equation*}
T_{1}
= \sqrt{\frac{2}{\xi}}
\sqrt{\frac{6 \beta}{\tau \lambda} + u(T)}
> \sqrt{\frac{2 u(T)}{\xi}}
\geq \sqrt{\frac{2 R}{\xi}}
= T.
\end{equation*}
Then 
we have
\begin{equation*}
u_{t}(t)
> \sqrt{\frac{\tau \lambda}{6}} \left( u(t) - u(T) \right)^{\frac{3}{2}}
\end{equation*}
and
\begin{equation*}
\left( u(t) - u(T) \right)^{-\frac{3}{2}} 
u_{t}(t) 
> \sqrt{\frac{\tau \lambda}{6}}
\end{equation*}
for all $t > T_{1}$.
Integrating this inequality
over $[T_{1}, t]$,
we have
\begin{eqnarray*}
0
< \frac{2}{\sqrt{u(t) - u(T) }}
< \frac{2}{\sqrt{u(T_{1}) - u(T) }}
- \sqrt{\frac{\tau \lambda}{6}} \left( t - T_{1} \right),
\end{eqnarray*}
which implies
that $\lim_{t \to T_{2}} u(t) = +\infty$,
where
\begin{eqnarray*}
T_{2}
= T_{1} + \sqrt{\frac{6}{\tau \lambda}} \frac{2}{\sqrt{u(T_{1}) - u(T) }}
< +\infty,
\end{eqnarray*}
contradicting the maximality of $T_{\infty} $ and hence 
necessarily $T_{\infty} < +\infty$.
Thanks to \eqref{eqn:integration-lyaunov},
we also have
\begin{equation}
\xi u + E_{\lambda}( u_{0}, v_{0} )
\leq - \int_{0}^{u} H_{\lambda}(s) \, ds + E_{\lambda}( u_{0}, v_{0} )
= \frac{1}{2} u_{t}^2,
\label{eqn:blowup-12}
\end{equation}
which implies 
that 
both $u$ and $u_{t}$ blow up to $+\infty$ as $t \to T_{\infty}$.

\medskip

\noindent {\bf III. The case $\alpha \neq 0$ and $I=(-\infty, b)$.}\\
If $v_{0} = 0$,
we have
\begin{equation}
u_{tt}(0) = -H_{\lambda}(u_{0}) \geq \xi > 0.
\label{eqn:blowup-1}
\end{equation}
Hence 
for $v_{0} \geq 0$,
we have $u_{t}(t)>0$
for sufficiently small $t>0$.
If there exists $T_{3} \in (0, T_{\infty})$
such that
$u_{t}(t)>0$ for all $t \in (0, T_{3})$ 
and
$u_{t}(T_{3})=0$,
then
we have
$u_{tt}(T_{3}) > 0$
similarly to \eqref{eqn:blowup-1},
which contradicts the positivity
of $u_{t}$.
Hence,
if necessary,
we retake the initial value
as $u_{0} = u(T_{4})$
and
$v_{0} = v(T_{4})$
for some $T_{4}>0$
so that
$u_{0} > 0$,
$v_{0} > 0$,
$u(t)>0$ 
and
$u_{t}(t)>0$
hold 
for all $t \in (0, T_{\infty})$.
We estimate $u_{t}(t)$.
First
if
$u_{tt}(0)>0$
holds,
we have
$u_{t}(t) \geq v_{0}$
for sufficiently small $t>0$.
On the other hand,
if $u_{tt}(t) \leq 0$
holds
for some $t \geq 0$,
we have
\begin{equation*}
\alpha \eta u_{t}^{\theta +1}
\geq u_{tt} + \alpha f(u_{t}) 
= -H_{\lambda}(u) 
\geq \xi 
> 0 
\end{equation*}
and then
\begin{equation*}
u_{t}
\geq \left( \frac{\xi}{\alpha \eta} \right)^{\frac{1}{\theta +1}}.
\end{equation*}
Thus 
we have
\begin{equation}
u_{t}
\geq \min \left\{ v_{0}, 
\left( \frac{\xi}{\alpha \eta} \right)^{\frac{1}{\theta +1}}
\right\}
\equiv C_{1}
> 0
\label{eqn:blowup-III-1}
\end{equation}
for all $t \in (0, T_{\infty})$,
which yields 
\begin{equation}
u(t)>u_{0}+C_{1}t
\label{eqn:blowup-III-2}
\end{equation}
for all $t \in (0, T_{\infty})$,
which brings back to the same situation
of case I above.

\medskip

\noindent {\bf IV. The case $\alpha \neq 0$ and $I={\mathbb R}$.}\\
By estimates carried out 
in the case III,
we have that 
\eqref{eqn:blowup-III-1}
and
\eqref{eqn:blowup-III-2} hold.
Hence,
$(u, u_{t})$ is unbounded
in ${\mathbb R}^2$
for $t \in (0, T_{\infty})$,
where
$T_{\infty} \leq +\infty$.
In the case
of $T_{\infty} < +\infty$,
both $u$ and $u_{t}$ blow up to $+\infty$ as $t \to T_{\infty}$
similarly to \eqref{eqn:blowup-12}.
Next let us consider the case $T_{\infty} = +\infty$ and let us prove
that
$u_{t}(t) \to +\infty$, as $ t \to +\infty$.
Now
suppose
that
there exists a constant $C_{2}>0$
satisfying 
\begin{equation*}
C_{1} < u_{t} < C_{2}
\end{equation*}
for all $t \geq T$.
Then, buying the line of  \eqref{eqn:blowup-10}
we obtain the following differential inequality
\begin{eqnarray*}
u_{tt} + \alpha \eta u_{t}^{\theta + 1} 
&\geq& u_{tt} + \alpha f( u_{t} ) \\
&>& \frac{\tau \lambda}{2} \left( u(t) - u(T) \right)^2 
- \beta \left( u(t) - u(T) \right) \\
&=& \frac{\tau \lambda}{2} 
\left( u(t) - u(T) - \frac{\beta}{\tau \lambda} \right)^2 
- \frac{\beta^2}{2 \tau \lambda}
\end{eqnarray*}
for all $t \geq T$.
Thus 
for sufficiently large $t>T$,
we have
\begin{equation*}
u_{tt} \geq \frac{\tau \lambda}{2} 
\left( C_{1}t +u_{0} - u(T) - \frac{\beta}{\tau \lambda} \right)^2 
- \frac{\beta^2}{2 \tau \lambda} - \alpha \eta C_{2}^{\theta + 1},
\end{equation*}
which yields $u_{t}(t) \to +\infty$,
as $ t \to +\infty$
by integration
and
contradicting the boundedness
of $u_{t}$.
Hence,
both $u$ and $u_{t}$ blow up to $+\infty$ as $t \to T_{\infty}$.
\hfill $\Box$

\vspace{1pc}

\noindent {\it Proof of Theorem \ref{thm:blowup-large-IV}.}
By hypothesis,
$\phi_{2} < u(t) < b$ holds for sufficiently small $t>0$.
If $v_{0}=0$,
we have
\begin{equation*}
u_{tt}(0) 
= -H_{\lambda}(u_{0}) 
> -H_{\lambda}(\phi_{2} )
= 0
\end{equation*}
by Corollary \ref{coro:H-sign} and Remark \ref{coro:H^*-sign}.
Hence
we may assume 
that
$\phi_{2} < u_{0} < u(t) < b$ 
and
$u_{t}(t) >0$
holds 
for all $t \in (0, T_{\infty})$.
Now set $\xi = - H_{\lambda}(u_{0}) > 0$ and proceed as in the proof of Theorem \ref{thm:blowup-large-para}.
\hfill $\Box$

\vspace{1pc}

\noindent Let us denote by $\gamma (t; \lambda, \alpha)$ 
the orbit of solution of \eqref{eqn:ODE-hyperbolic}
starting at $(0, 0)$.
We note that from Lemma \ref{lem:eqilibria-sink-centre}, one has 
 $\alpha f'(0) = 0$ 
if and only if 
$( \phi_{1}(\lambda), 0)$ is a centre.

\begin{proposition}
Let $\alpha = 0$. We have that 
$0 < \lambda < {\overline \lambda}$ is satisfied
if and only if the orbit $\gamma (t; \lambda, 0)$ is periodic.
If ${\overline \lambda} < \lambda < \lambda^*$,
then the orbit is unbounded.
\label{prop:ex-periodic-sol}
\end{proposition}
{\it Proof.}
In the conservative case $\alpha = 0$
we have
\begin{equation*}
E_{\lambda}( u, u_{t} )
= E_{\lambda}( 0 ,0 )
\Longleftrightarrow 
J_{\lambda}(u) + \frac{1}{2}u_{t}^2
= 0
\label{eqn:sol-est-3}
\end{equation*}
by \eqref{eqn:integration-lyaunov}.
Since $\left( u_{t} \right)_{t} = -  H_{\lambda}(u)$ is positive either for $u < \phi_{1}$ or $\phi_{2} < u < b$ and negative
for $\phi_{1} < u < \phi_{2}$, the orbit is periodic
if and only if
there exists $u_{1} \in (\phi_{1}, \phi_{2})$
such that it crosses the $u$-axis at the point $\left( u_{1},0 \right)$
by Remark \ref{rem:ex-global} 
and 
Theorem \ref{thm:blowup-large-IV},
which implies that
\begin{equation*}
J_{\lambda}( u_{1} )=0.
\end{equation*}
This condition is equivalent to $\lambda \in ( 0, {\overline \lambda} )$
by the monotonicity of $J_{\lambda}(u)$ with respect to $u$
in Corollary \ref{coro:H-sign},
Lemmas \ref{lem:local-min-max}
and 
\ref{lem:local-max}.
If ${\overline \lambda} < \lambda < \lambda^*$,
there does not exist such $u_{1}$,
which means that $\phi_{2} < u < b$ and $u_{t} > 0$ for sufficiently large $t$.
Finally apply Theorem \ref{thm:blowup-large-IV} to conclude.
\hfill $\Box$

\vspace{1pc}

\begin{remark}
Since $\lambda < {\overline \lambda}$ is equivalent to $l(\lambda)<0$,
we have $\left( 0, 0 \right) \in D_{0}$.
Thus 
by Theorem \ref{thm:ex-global-sol},
the orbit $\gamma (t; \lambda, \alpha)$ for $\alpha > 0$ exists globally in time 
and 
converges to $(\phi_{1}(\lambda), 0)$, as $t \to +\infty$.
Hence along with Proposition \ref{prop:ex-periodic-sol},
we have ${\overline \lambda} \leq \lambda \left( 0, 0 \right) (\alpha)$ 
for $\alpha > 0$ 
and 
${\overline \lambda} = \lambda \left( 0, 0 \right) (0) < \lambda^*$.
\label{rem:ex-global-sol-origin}
\end{remark}

\section{Properties of the dissipative orbit}
\label{sec:properties of the orbit}

In this Section,
we study the properties of the orbit starting at the origin 
for the dissipative case,
that is,
$\alpha > 0$.
Moreover,
we also assume
that
$f(v)=v$.
Then,
$\left( \phi_{1} (\lambda), 0 \right)$
is a hyperbolic sink
for all $\lambda \in ( {\overline \lambda},  \lambda^* )$.
The argument proceeds in the same way as in Section 3 of \cite{MR3592652}.

\vspace{1pc}

Let us first define a few sets which will be used in the sequel:
\begin{align*}
&\Gamma :=\left\{ \left( \lambda, \alpha \right) \in {\mathbb R}^2 \mid
{\overline \lambda} < \lambda < \lambda^*, \ \alpha > 0 \right\}; \\
& \Gamma_{1} := \left\{ \left( \lambda, \alpha \right) \in \Gamma \mid
\gamma (t; \lambda, \alpha) \to \left( \phi_{1} (\lambda), 0 \right),
\ \text{as} \ t \to +\infty \right\};\\
& \Gamma_{2} := \left\{ \left( \lambda, \alpha \right) \in \Gamma \mid
\gamma (t; \lambda, \alpha) \to \left( \phi_{2} (\lambda), 0 \right),
\ \text{as} \ t \to +\infty
\right\}; \\
& \Gamma_{3} := \left\{ \left( \lambda, \alpha \right) \in \Gamma \mid
\gamma (t; \lambda, \alpha) \ \text{becomes unbounded}
\ \text{as} \ t \to T_{\infty}
\right\},
\end{align*}
where $\gamma (t; \lambda, \alpha)$ is the orbit of the solution $(u, v) = (u, u_{t})$ 
of \eqref{eqn:ODE-hyperbolic} starting at $(0, 0)$.
In what follows we will also use for convenience the following equivalent notations 
\begin{equation*}
\gamma (t)
=\gamma (t; \lambda, \alpha)
= \left[ \begin{array}{c}
u (t; \lambda, \alpha)\\
v (t; \lambda, \alpha)
\end{array} \right]
= \left[ \begin{array}{c}
u (t)\\
v (t)
\end{array} \right]
\label{eqn:orbit-alpha}
\end{equation*}
for any fixed $\alpha,\lambda >0$.

\begin{lemma}
$\Gamma_{1}$ and $\Gamma_{3}$ are open subsets of $\Gamma$.
\label{lem:openness-parameter-region}
\end{lemma}
{\it Proof.}
If we take $\left( \lambda_{0}, \alpha_{0} \right) \in \Gamma_{1}$,
then
$\left( u, v \right) = \left( \phi_{1} (\lambda_{0}), 0 \right)$ is a hyperbolic sink 
by Lemma \ref{lem:eqilibria-sink-centre}.
Hence,
there exists $r>0$ 
depending only on $\lambda_{0}$, $\alpha_{0}$, $\beta$ and $g$ 
such that
we can find an invariant disk $B_{r} \left( \phi_{1} (\lambda_{0}), 0 \right)$ 
of radius $r$ centred 
at $\left( \phi_{1} (\lambda_{0}), 0 \right)$
for the dynamical system 
induced by \eqref{eqn:dyn-sys}.
In other words,
there exists $T>0$ such that
$\gamma (t; \lambda_{0}, \alpha_{0}) \subset B_{r} \left( \phi_{1} (\lambda_{0}), 0 \right)$ 
for all $t>T$.
If necessary,
we can take $r>0$ so small such that 
\begin{equation*}
B_{r} \left( \phi_{1} (\lambda_{0}), 0 \right)
\subset \left\{ \left( u,v \right) \in {\mathbb R}^2 \mid
u \leq p \right\},
\end{equation*}
where $p$ is defined in Theorem \ref{thm:stationary-sol} and
satisfies $\phi_{1} (\lambda_{0}) < p < \phi_{2} (\lambda_{0})$.
For $\left( \lambda, \alpha \right) \in \Gamma$ 
sufficiently close to $\left( \lambda_{0}, \alpha_{0} \right)$,
$\gamma (t; \lambda, \alpha) \subset B_{r} \left( \phi_{1} (\lambda_{0}), 0 \right)$ 
for sufficiently large $t$
thanks to the continuous dependence on parameters.
Since $\left( \phi_{1} (\lambda), 0 \right)$ is also a hyperbolic sink,
$\gamma (t; \lambda, \alpha)$ converges $\left( \phi_{1} (\lambda), 0 \right)$ 
as $t \to +\infty$.
Hence
we have $\left( \lambda, \alpha \right) \in \Gamma_{1}$.\\
Now take $\left( \lambda_{0}, \alpha_{0} \right) \in \Gamma_{3}$ and define
\begin{equation*}
{\mathcal U}_{\lambda, \alpha} := \left\{ \left( u,v \right) \in {\mathbb R}^2
\mid \phi_{2}(\lambda) < u < b, \
0 < v < -\frac{1}{\alpha} H_{\lambda}(u) \right\}.
\end{equation*}
For sufficiently large $t>0$,
we have $\gamma (t; \lambda_{0}, \alpha_{0}) \subset {\mathcal U}_{\lambda_{0}, \alpha_{0}}$
by a phase plane analysis of \eqref{eqn:dyn-sys}
togehter with $u_{tt}(0) = \lambda_{0} > 0$ and Remark \ref{rem:ex-global}.
Again, by continuous dependence,
for large $t>0$
we have $\gamma (t; \lambda, \alpha) \subset {\mathcal U}_{\lambda_{0}, \alpha_{0}}$
for $\left( \lambda, \alpha \right) \in \Gamma$ 
sufficiently close to $\left( \lambda_{0}, \alpha_{0} \right)$. Since we have $\phi_{2} (\lambda) < u (t; \lambda, \alpha) < b$ 
and 
$u_{t} (t; \lambda, \alpha) > 0$ 
for sufficiently large $t$,
$\gamma (t; \lambda, \alpha)$ becomes unbounded
by Theorem \ref{thm:blowup-large-IV},
which proves that 
$(\lambda, \alpha) \in \Gamma_{3}$.
\hfill $\Box$

\vspace{1pc}

\begin{proposition}\label{prop:parameter-rigion-decomopsition}
We have $\Gamma
= \bigcup_{i=1}^3 \Gamma_{i}$.
\end{proposition}
{\it Proof.}
Note 
that 
every bounded orbit for $\alpha > 0$ converges to $\phi_{1}$ or $\phi_{2}$ 
as $t \to +\infty$
by Corollary 8.5.1 
in \cite{MR3380967}.
If an unbounded orbit $\gamma (t; \lambda, \alpha)$ exists,
then
$u(t) \geq 0$ and $u_{t}(t) \geq 0$ hold 
for sufficiently small $t \geq 0$ by $u_{tt}(0) = \lambda > 0$
and
the orbit $\gamma(t; \lambda, \alpha)$ enters ${\mathcal U}_{\lambda, \alpha}$
by Remark \ref{rem:ex-global}.
Then
we have $(\lambda, \alpha) \in \Gamma_{3}$ 
as in the proof of Theorem \ref{thm:blowup-large-IV}.
\hfill $\Box$

\vspace{1pc}

\noindent 
Let $m$ be a negative constant to be determined later.
Let us define the line segments $s_{i}$ for $i=1$,$2$,$3$ and the triangular region $\mathcal{T}$ as follows:
\begin{align*}
& s_{1} :=\left\{ \left( u, v \right) = \left( \phi_{2} \left( \lambda \right), v \right)
\mid 
m \left( \phi_{2} \left( \lambda \right) - \phi_{1} \left( \lambda \right) \right) 
\leq v \leq 0 \right\}; \\
& s_{2} := \left\{ \left( u, v \right) = \left( u, 0 \right) 
\mid
\phi_{1} \left( \lambda \right) \leq u \leq \phi_{2} \left( \lambda \right) \right\};\\
& s_{3} := \left\{ \left( u, v \right) 
= \left( u, m \left( u - \phi_{1} \left( \lambda \right) \right) \right) 
\mid
\phi_{1} \left( \lambda \right) \leq u \leq \phi_{2} \left( \lambda \right) \right\};\\
& \mathcal{T} := \left\{ \left( u, v \right) \in {\mathbb R}^2 
\mid
\phi_{1} \left( \lambda \right) \leq u \leq \phi_{2} \left( \lambda \right), \ 
m \left( u - \phi_{1} \left( \lambda \right) \right) \leq v \leq 0 \right\} \ .
\end{align*}

\begin{proposition}
For $\alpha > 2 \sqrt{H_{\lambda}'(\phi_{1})}$,
there exists a heteroclinic orbit 
from $\left( \phi_{2}, 0 \right)$ to $\left( \phi_{1}, 0 \right)$.
\label{prop:ex-heteroclinic-orbit}
\end{proposition}
{\it Proof.}
Since the vector field on $s_{1}$ and $s_{2}$ 
defined by \eqref{eqn:dyn-sys} points inward $\mathcal{T}$,
we can choose $m<0$ such that the vector field on $s_{3}$ also points inward $\mathcal{T}$.
Denote by $\bm{N}$ the normal vector on $s_{3}$ 
\begin{equation*}
\bm{N}
= \left[ \begin{array}{cc}
m\\
-1
\end{array} \right]
\label{eqn:normal-l-3}
\end{equation*}
and by $\bm{V}$ the vector on $s_{3}$ defined as follows
\begin{equation*}
\bm{V}
= \left[ \begin{array}{cc}
u_{t}\\
-\alpha u_{t} - H_{\lambda}(u)
\end{array} \right]
= \left[ \begin{array}{cc}
m \left( u - \phi_{1} \right) \\
- \alpha m \left( u - \phi_{1} \right) - H_{\lambda}(u)\\
\end{array} \right]
\label{eqn:vector-field-l-3}
\end{equation*}
for $\phi_{1} < u < \phi_{2}$.
Let us compute the inner product between $\bm{N}$ and $\bm{V}$ to obtain
\begin{eqnarray*}
\bm{N} \cdot \bm{V}
&=& m^2 \left( u - \phi_{1} \right) 
+ \alpha m \left( u - \phi_{1} \right) 
+ H_{\lambda}(u) \\
&=& \left( u - \phi_{1} \right) 
\left( m^2 + \alpha m + \frac{H_{\lambda}(u)}{ u - \phi_{1}} \right) \\
&\leq& \left( u - \phi_{1} \right) \left( m^2 + \alpha m + H_{\lambda}'(\phi_{1}) \right)
\label{eqn:inner-product}
\end{eqnarray*}
by Lemma \ref{eqn:condition-heteroclinic}.
Take $m$ such that $\bm{N} \cdot \bm{V} < 0$, 
that is, 
$\mu_{1}^{-} < m < \mu_{1}^{+}$,
where $\mu_{1}^{\pm}$ are defined in Section \ref{sec:dynamical system}
with
$f'(0)=1$.
Now the branch of the unstable manifold of $\left( \phi_{2}, 0 \right)$ that points inside the region
\begin{equation*}
\left\{ \left( u, v \right) \in {\mathbb R}^2
\mid u < \phi_{2} \left( \lambda \right), \ v < 0 \right\}
\label{eqn:invariant-region}
\end{equation*}
enters $\mathcal{T}$ and does not leave it.
Hence 
this bounded orbit in $\mathcal{T}$ converges to $\left( \phi_{1}, 0 \right)$ 
as $t \to +\infty$.
\hfill $\Box$

\vspace{1pc}

\begin{remark}
A typical example for $f$ is given by $f(v)=\left\vert v \right\vert^{\gamma}$,
where $\gamma  \geq 1$.
However,
we consider the case
of $f(v)=v$.
Indeed, if we take $f(v)=\left\vert v \right\vert^{2}$,
Proposition \ref{prop:ex-heteroclinic-orbit}
does not hold
because
\begin{eqnarray*}
\bm{N} \cdot \bm{V}
&=& m^2 \left( u - \phi_{1} \right) 
+ \alpha m^2 \left( u - \phi_{1} \right)^2 
+ H_{\lambda}(u)
> 0
\label{eqn:inner-product}
\end{eqnarray*}
for any $m \in {\mathbb R}$
and
$u \in (\phi_{1}, \phi_{2})$.
\label{rem:choice-g}
\end{remark}

\vspace{1pc}

\noindent In order to prove the monotonicity of the orbit,
let us introduce the following notation 
\begin{equation*}
\gamma_{i} (t)
=\gamma (t; \lambda_{i}, \alpha)
= \left[ \begin{array}{cc}
u (t; \lambda_{i}, \alpha)\\
v (t; \lambda_{i}, \alpha)
\end{array} \right]
= \left[ \begin{array}{cc}
u_{i} (t)\\
v_{i} (t)
\end{array} \right]
\label{eqn:orbit-i-alpha}
\end{equation*}
for any fixed $\alpha >0$ and $0 < \lambda_{i} < \lambda^*$, $i=1$, $2$.
We next prove the monotonicity of the orbit in $\lambda$ for fixed $\alpha$.

\begin{proposition}
For fixed $\alpha >0$,
$v (t; \lambda, \alpha)$ is increasing 
with respect to  $\lambda$
as long as $v (t; \lambda, \alpha)>0$.
Moreover,
if there exists $\lambda_{0} \in (0, \lambda^*)$ 
such that $(\lambda_{0}, \alpha) \in \Gamma_{2} \cup \Gamma_{3}$,
then
we have $(\lambda, \alpha) \in \Gamma_{3}$ for all $\lambda > \lambda_{0}$.
\label{prop:monotonicity-orbit-alpha}
\end{proposition}
{\it Proof.}
Since we have
$v_{i} (0) = 0$ 
and 
$(v_{i})_{t} (0) = \lambda_{i}$,
$v_{1} (t) < v_{2} (t)$ holds
for sufficiently small $t>0$.
For the second component
of vector field 
defined in \eqref{eqn:dyn-sys},
we have
\begin{equation*}
\frac{d}{d \lambda}
\left( - \alpha v - H_{\lambda} ( u ) \right)
= g(u)
>0.
\end{equation*}
Let us prove the second statement.
If $(\lambda_{0}, \alpha) \in \Gamma_{3}$,
the monotonicity yields $(\lambda, \alpha) \in \Gamma_{3}$ for all $\lambda > \lambda_{0}$.
If $(\lambda_{0}, \alpha) \in \Gamma_{2}$,
the monotonicity
of $\phi_{i}(\lambda)$ and that of the orbit 
in $\lambda$
imply
$(\lambda, \alpha) \in \Gamma_{3}$ 
by Theorem \ref{thm:blowup-large-IV}.
\hfill $\Box$

\vspace{1pc}

\noindent As stated after the proof of Proposition 6 in \cite{MR3592652},
the stable local manifold of the saddle $\left( \phi_{2}(\lambda), 0 \right)$ plays 
a crucial role in determining the threshold of the parameter $\lambda$.
Now for fixed $\lambda \in ({\overline \lambda}, \lambda^*)$,
we regard the behaviour of the local stable manifold as a function of $\alpha$.
We shall prove that the manifold crosses the positive $u$-axis for small $\alpha>0$.
Then the orbit $\gamma (t)$ cannot approach the stationary points,
which implies that $(\lambda, \alpha) \in \Gamma_{3}$.
On the other hand,
for large $\alpha>0$,
we prove that the manifold crosses the negative $u$-axis.
In this case,
the solution $(u, v)$ is bounded 
for all $t \geq 0$ 
and 
$(\lambda, \alpha) \not \in \Gamma_{2}$,
which means that $(\lambda, \alpha) \in \Gamma_{1}$.
Finally,
we uniquely determine $\alpha^*(\lambda) >0$ 
such that 
the manifold crosses the $u$-axis at $u=0$.
Then we establish the continuity and monotonicity of $\alpha^*(\lambda) >0$ 
with respect to $\lambda$
and define $\lambda^* \left( 0, 0 \right) (\alpha)$ 
by the inverse function of $\alpha^*(\lambda)$.
In order to analyze the stable manifold of $\left( \phi_{2}(\lambda), 0 \right)$,
we perform the following change of variables in \eqref{eqn:ODE-hyperbolic} 
\begin{equation*}
\left\{\begin{array}{ll}
t=-s, \\
U(s) = \phi_{2}(\lambda) - u(t), \\
V(s) = v(t)
\end{array}\right.
\label{eqn:transform}
\end{equation*}
and consider
\begin{equation}
\frac{d}{ds}
\left[ \begin{array}{cc}
U\\
V
\end{array} \right]
= \left[ \begin{array}{cc}
V\\
\alpha V + H_{\lambda} \left( \phi_{2}(\lambda) - U \right) 
\end{array} \right]
\label{eqn:transformed-ODE-hyperbolic}
\end{equation}
for $s<0$ 
with $U(-\infty) = V(-\infty) = 0$.
As we have seen in Section \ref{sec:dynamical system},
the eigenvalues $\eta^{\pm}$ of the linearized operator at $(U,V)=(0,0)$ 
corresponding to \eqref{eqn:transformed-ODE-hyperbolic} are given by
\begin{equation*}
\eta^{\pm}
= \eta^{\pm}(\alpha)
= \frac{\alpha \pm \sqrt{\alpha^2 - 4 H_{\lambda}'(\phi_{2})}}{2}.
\label{eqn:trnsformed-liearized-op}
\end{equation*}
The branch of the local unstable manifold can be expressed by the graph 
$V = \Phi (U; \lambda, \alpha)$
as long as $V(s)>0$
since $V(s)=U_{s}(s)>0$.
First,
it is clear that $\Phi (0; \lambda, \alpha)=0$.
Since we concentrate on the unstable manifold,
we have $\left( d \Phi / dU \right) (0; \lambda, \alpha)=\eta^+$
by $H_{\lambda}'(\phi_{2}) < 0$.
Finally, we have
\begin{equation}
\frac{d \Phi}{dU}(U; \lambda, \alpha)
= \frac{V_{s}}{U_{s}}
= \alpha + \frac{H_{\lambda} \left( \phi_{2}(\lambda) - U \right)}{\Phi (U; \lambda, \alpha)}.
\label{eqn:derivative-unstable-mfd}
\end{equation}
We denote
$\Phi_{i} (U)
= \Phi (U; \lambda, \alpha_{i})$
for fixed $\lambda \in ({\overline \lambda}, \lambda^*)$ 
and 
$\alpha_{i}>0$ 
with $i=1$, $2$
and establish the monotonicity of $\Phi (U; \lambda, \alpha)$ 
with respect to $\alpha$.

\begin{proposition}
For fixed $\lambda \in ( 0, \lambda^*)$,
let $0 < \alpha_{1} < \alpha_{2}$.
The graph $V = \Phi_{2}(U)$ stays above the graph $V =  \Phi_{1}(U)$.
Moreover,
the graph $V = \Phi_{2}(U)$ and $V =  \Phi_{1}(U)$ never intersect each other
as long as they are defined.
\label{prop:monotonicity-unstable-mfd-lambda}
\end{proposition}
{\it Proof.}
Since we have $\Phi_{2} (0) - \Phi_{1} (0) = 0$ and 
$\left( d \Phi_{2} / dU \right) (0) - \left( d \Phi_{1} / dU \right) (0) 
= \eta^+(\alpha_{2}) - \eta^+(\alpha_{1}) > 0$,
we obtain
$\Phi_{2}(U) - \Phi_{1}(U) > 0$ for sufficiently small $U>0$.
Thanks to \eqref{eqn:derivative-unstable-mfd},
$V = \Phi_{2}(U)$ and $V =  \Phi_{1}(U)$ can not intersect each other.
\hfill $\Box$

\medskip

\noindent For fixed $\lambda \in ({\overline \lambda}, \lambda^*)$,
let
\begin{equation*}
K(\lambda):= \left\{ \alpha \geq 0
\mid \ 
\text{there exists} \ P_{\alpha}>0 \ \text{such that} \ 
\Phi \left( P_{\alpha}; \lambda, \alpha \right) = 0 \right\}.
\label{eqn:set-I}
\end{equation*}
We are interested in the set $L(\lambda)$ of all points $P_{\alpha}$ 
defined by
\begin{equation*}
L(\lambda):= \left\{ P_{\alpha}
\mid \alpha \in K(\lambda) \right\}.
\label{eqn:set-J}
\end{equation*}
$L(\lambda)$ consists of the points 
where the unstable manifold intersects the positive $U$-axis.
We shall show that $K(\lambda)$ is a non-empty interval 
and that 
$L(\lambda)$ is an unbounded interval.
For this purpose,
let us define two lines parallel to the $V$-axis as follows:
\begin{align*}
& M(\lambda):= \left\{ \left( \phi_{2} (\lambda) - \phi_{1} (\lambda), V \right) 
\mid V \geq 0 \right\}; \\
& M_{\Phi}(\lambda):= \left\{ \left( \phi_{2} (\lambda) - \phi_{1} (\lambda), 
\Phi \left( \phi_{2} (\lambda) - \phi_{1} (\lambda); \lambda, \alpha \right) \right) 
\mid \alpha \geq 0 \right\}. 
\label{eqn:lines-K-L}
\end{align*}
Thanks to the continuity of the intersection point 
with respect to $\alpha$ proved in \cite{MR511133},
$M_{\Phi}(\lambda)$ is an interval in $M(\lambda)$.

\begin{lemma}\label{lem:structure-K}
$M_{\Phi}(\lambda)
= \left\{ \phi_{2} (\lambda) - \phi_{1} (\lambda) \right\} \times
[ \Phi \left( \phi_{2} (\lambda) - \phi_{1} (\lambda); \lambda, 0 \right), +\infty)$ 
for $\lambda \in ({\overline \lambda}, \lambda^*)$.
\end{lemma}
{\it Proof.}
Noting that 
$H_{\lambda} \left( \phi_{2} - U \right) \geq 0$ for $U \in [0, \phi_{2} - \phi_{1}]$,
we have
\begin{equation*}
\frac{d \Phi}{dU}(U; \lambda, \alpha)
\geq \alpha 
\label{eqn:structure-K-ODE}
\end{equation*}
for any $\alpha \geq 0$ and $U \in [0, \phi_{2} - \phi_{1}]$
by \eqref{eqn:derivative-unstable-mfd}.
Hence integrating this inequality over $[0, \phi_{2} - \phi_{1}]$, 
we obtain
\begin{equation*}
\Phi \left( \phi_{2} - \phi_{1}; \lambda, \alpha \right)
\geq \alpha \left( \phi_{2} - \phi_{1} \right)
\to + \infty, \text{ as } \alpha \to +\infty \ .
\label{eqn:structure-K-est-V}
\end{equation*}
\hfill $\Box$

\vspace{1pc}

\noindent Next we establish a few properties of $K(\lambda)$ and $L(\lambda)$, for which we set $A(\lambda) := \sup K(\lambda)$ and $P(\lambda) := \sup L(\lambda)$.

\begin{lemma}
$K(\lambda)$ and $L(\lambda)$ are nonempty intervals
for $\lambda \in ({\overline \lambda}, \lambda^*)$.
\label{lem:interval-I-J}
\end{lemma}
{\it Proof.}
There exists some $u \in (0, \phi_{1})$
such that
$J_{\lambda}(u) = J_{\lambda}(\phi_{2})$
holds.
From the arguments of phase plane analysis 
in Theorem \ref{thm:ex-global-sol},
Remark \ref{rem:ex-global}
and
Proposition \ref{prop:ex-periodic-sol},
$0 \in K(\lambda)$ holds.
Hence
$K(\lambda) \neq \emptyset$ follows
for $\lambda \in ({\overline \lambda}, \lambda^*)$.
Finally,
as mentioned in Proposition \ref{prop:monotonicity-unstable-mfd-lambda},
the continuity and monotonicity of $\Phi$ yield the conclusion.
\hfill $\Box$

\vspace{1pc}

\noindent We see that $K(\lambda)$ and $L(\lambda)$ are intervals.
Next we prove that the right endpoint of $K(\lambda)$ is open.

\begin{lemma}
$K(\lambda) = [0, A(\lambda) )$
for $\lambda \in ({\overline \lambda}, \lambda^*)$.
\label{lem:structure-I}
\end{lemma}
{\it Proof.}
We may assume $A(\lambda) < +\infty$. Suppose by contradiction that $K(\lambda) = [0, A(\lambda) ]$.
The orbit $\Phi \left( U; \lambda, A(\lambda) \right)$ intersects the line $M(\lambda)$.
Then,
by definition
we can find $P_{A(\lambda)} > 0$ 
such that 
$\Phi \left( P_{A(\lambda)}; \lambda, A(\lambda) \right) = 0$.
Next, 
the orbit starting at $(P_{A(\lambda)},0)$ enters the region 
$\left\{ \left( U, V \right) \in {\mathbb R}^2
\mid V<0 \right\}$
with the $U$-coordinate decreasing.
Since the set
\begin{equation*}
\left\{ \left( \phi_{2} (\lambda) - \phi_{1} (\lambda), 
\Phi \left( \phi_{2} (\lambda) - \phi_{1} (\lambda); \lambda, \alpha \right) \right)
\mid \alpha \in [0, A(\lambda) + \varepsilon ] \right\}
\label{eqn:def-subset-K}
\end{equation*}
for any $\alpha \geq 0$ 
and
sufficiently small $\varepsilon >0$
describes a closed interval in the line $M(\lambda)$
by the continuity and monotonicity of $\Phi$ with respect to $\alpha$ 
and 
Lemma \ref{lem:structure-K},
the behaviour of the orbit $\Phi \left( U; \lambda, \alpha \right)$ is 
the same as that of $\Phi \left( U; \lambda, A(\lambda) \right)$
for $\alpha \in ( A(\lambda), A(\lambda) + \varepsilon)$ 
with sufficiently small $\varepsilon > 0$ 
by the continuous dependence 
on the parameter
and
initial value.
Hence
we obtain $P_{\alpha}$ 
satisfying $\Phi \left( P_{\alpha}; \lambda, \alpha \right) = 0$
for $\alpha \in ( A(\lambda), A(\lambda) + \varepsilon)$,
which contradicts the definition of $A(\lambda)$.
\hfill $\Box$

\vspace{1pc}

\noindent Finally,
we prove that $L(\lambda)$ is an unbounded interval.

\begin{lemma}
$P(\lambda) = + \infty$
for $\lambda \in ({\overline \lambda}, \lambda^*)$.
\label{lem:structure-J}
\end{lemma}
{\it Proof.}
First,
we deal with the case of $A(\lambda) < +\infty$.
Assume by contradiction that $P(\lambda) < + \infty$.
The orbit $\Phi \left( U; \lambda, \alpha \right)$ intersects the line $M(\lambda)$ 
for any $\alpha \in ( 0, A(\lambda) )$.
In addition,
due to \eqref{eqn:derivative-unstable-mfd},
$\left( d \Phi \right) / \left( dU \right) (U; \lambda, \alpha) < +\infty$ 
on every finite interval 
for any $\alpha \in ( 0, A(\lambda) )$
as long as $\Phi (U; \lambda, \alpha) > 0$.
Note that
$A(\lambda) \not \in K(\lambda)$.
Thus
the continuous dependence on the parameter yields
$\Phi (U; \lambda, A(\lambda)) > 0$
for all $U > 0$.
Because of the continuous dependence,
for any $\alpha \in ( A(\lambda)  - \varepsilon, A(\lambda) )$ 
with sufficiently small $\varepsilon>0$,
we have
$\Phi (U; \lambda, \alpha) > 0$
for all $U \geq 0$,
which contradicts $\alpha \in K(\lambda)$.\\
Next,
we treat the case of $A(\lambda) = +\infty$.
For all $\alpha \geq 0$,
there exists $P_{\alpha} \in (\phi_{2} - \phi_{1}, +\infty)$ such that
we have $\Phi (P_{\alpha}; \lambda, \alpha) = 0$.
The intersection point of $V = \Phi (U; \lambda, \alpha)$ 
and 
$M(\lambda)$ lies in the region $\left\{ \left( U, V \right) \mid V' > 0 \right\}$.
The graph of $V = \Phi (U; \lambda, \alpha)$ must leave the region 
at the point $\left( U_{\alpha}, V_{\alpha} \right)$ 
satisfying
\begin{equation*}
V_{\alpha} = -\frac{1}{\alpha}H_{\lambda} \left( \phi_{2} - U_{\alpha} \right),
\quad
\phi_{2} - \phi_{1} < U_{\alpha} < P_{\alpha}
\quad \text{and} \quad 
\Phi ( \phi_{2} - \phi_{1} ; \lambda, \alpha) < V_{\alpha}.
\label{eqn:intersection-pt1}
\end{equation*}
Then
we have
\begin{equation*}
H_{\lambda} \left( \phi_{2} - U_{\alpha} \right)
< -\alpha \Phi ( \phi_{2}  - \phi_{1}; \lambda, \alpha) \to -\infty
\label{eqn:intersection-pt2}
\end{equation*}
as $\alpha \to +\infty$
by Lemma \ref{lem:structure-K},
which implies that
$b > p > \phi_{1} > \phi_{2} - U_{\alpha} \to -\infty$
as $\alpha \to +\infty$,
where $p$ is defined in Theorem \ref{thm:stationary-sol}.
Eventually
we obtain
$P_{\alpha} > U_{\alpha} \to \infty$,
as $\alpha \to +\infty$.
\hfill $\Box$

\section{Proof of Theorem \ref{thm:dynamical-threshold}}
\label{sec:proof-of-main-TH}

The strategy to prove our main result is the following. 
For fixed $\lambda \in ( \overline{\lambda}, \lambda^*)$,
there exists $\alpha^*(\lambda)$ such that
the solution of \eqref{eqn:ODE-hyperbolic} exists globally in time 
for $\alpha^*(\lambda) < \alpha < +\infty$
and that is unbounded  
for $0 < \alpha < \alpha^*(\lambda)$.
Then,
the function $\alpha^*(\lambda)$ turns out to be  monotone increasing and continuous 
with respect to $\lambda$.
Finally,
we have $\alpha^*(\lambda) \to +\infty$, as $\lambda \to \lambda^*$.
We denote the inverse function of $\alpha^*(\lambda)$ 
by $\lambda \left( 0,0 \right) (\alpha)$ which we prove to have the desired properties as in Theorem  \ref{thm:dynamical-threshold}.

\vspace{1pc}

\begin{lemma}
For fixed $\lambda \in ( \overline{\lambda}, \lambda^* )$,
there exists $\alpha^*(\lambda) > 0$ such that
the solution of \eqref{eqn:ODE-hyperbolic} exists globally in time 
for $\alpha^*(\lambda) < \alpha < +\infty$
and that 
becomes unbounded 
for $0 < \alpha < \alpha^*(\lambda)$.
\label{lem:ex-thresh-ld-alpha}
\end{lemma}
{\it Proof.}
For any $\alpha \in (0, A(\lambda) )$,
we have $P_{\alpha} > P_{0}$ such that $\Phi \left( P_{\alpha}; \lambda, \alpha \right) = 0$ 
by Proposition \ref{prop:monotonicity-unstable-mfd-lambda}.
Let 
\begin{equation*}
x_{\alpha} = \phi_{2} - P_{\alpha}.
\label{eqn:def-x-alpha}
\end{equation*}
In other words,
$x_{\alpha}$ is an intersection point between the stable manifold 
of the saddle $( \phi_{2}(\lambda), 0)$ 
and 
the $u$-axis.
Since the point $P_{0}$ is determined by the orbit of the conservative case,
we have $x_{0}>0$,
that is,
$P_{0} < \phi_{2}$
by Proposition \ref{prop:ex-periodic-sol} 
for the case  $\lambda \in ( \overline{\lambda}, \lambda^* )$.
By Lemma \ref{lem:structure-J}, 
the monotonicity and continuity of the intersection point,
there exists $\alpha^*(\lambda) > 0$ uniquely determined 
such that 
$x_{\alpha^*(\lambda)}=0$, or $P_{\alpha^*(\lambda)} = \phi_{2}$.
First,
we consider the case of $0 < \alpha < \alpha^*(\lambda)$.
Since we have
\begin{equation*}
\phi_{2} - \phi_{1}
< P_{0}
< P_{\alpha}
< \phi_{2}
\Longleftrightarrow
0
< x_{\alpha}
< x_{0}
< \phi_{1},
\label{eqn:intersection-point-1}
\end{equation*}
$(0,0)$ is not in the domain of attraction of $( \phi_{1}(\lambda), 0)$ 
for the orbit $\gamma (t; \lambda, \alpha)$.
Hence
it enters
\begin{equation*}
Z_{1}:= \left\{ \left( u, v \right) \in {\mathbb R}^2
\mid
0 < u < b, \ v > 0, \ v_{t} > 0 \right\}
\label{eqn:blow-up-region1}
\end{equation*}
and
\begin{equation*}
Z_{2}:=\left\{ \left( u, v \right) \in {\mathbb R}^2
\mid
0 < u < b, \ v > 0, \ v_{t} < 0 \right\}.
\label{eqn:blow-up-region2}
\end{equation*}
Then,
we can find $T>0$ such that $\phi_{2}(\lambda) < u(T) < b$, $v(T)>0$ and $v_{t}(T)=0$
by Proposition \ref{prop:monotonicity-unstable-mfd-lambda}.
Hence
the solution blows up by Theorem \ref{thm:blowup-large-IV}.

\medskip

\noindent Next,
we consider the case $\alpha > \alpha^*(\lambda)$.
Similarly,
we obtain 
\begin{equation*}
\phi_{2}
< P_{\alpha}
\Longleftrightarrow
x_{\alpha}
< 0.
\label{eqn:stable-case1}
\end{equation*}
Threfore, $\gamma (t; \lambda, \alpha)$ enters $Z_{1}$ 
and 
then $Z_{2}$.
Since $\gamma (t; \lambda, \alpha)$ stays 
below the stable manifold of $( \phi_{2}(\lambda), 0)$,
it necessarily intersects the $u$-axis 
between $\phi_{1}(\lambda)$ 
and 
$\phi_{2}(\lambda)$.
Hence
the solution exists globally in time by Remark \ref{rem:ex-global}.
\hfill $\Box$

\vspace{1pc}

\begin{lemma}
$\alpha^*(\lambda)$ 
is monotone increasing
for $\lambda \in ( \overline{\lambda}, \lambda^*)$.
\label{lem:monotonicity-thresh-ld-alpha}
\end{lemma}
{\it Proof.}
Let $\overline{\lambda} < \lambda_{1} < \lambda_{2} < \lambda^*$.
We have $(\lambda_{1}, \alpha^* (\lambda_{1}) ) \in \Gamma_{2}$
by Proposition \ref{prop:parameter-rigion-decomopsition}
and moreover $(\lambda_{2}, \alpha^* (\lambda_{1}) ) \in \Gamma_{3}$
by Proposition \ref{prop:monotonicity-orbit-alpha}.
Hence,
from Lemma \ref{lem:ex-thresh-ld-alpha},
$\alpha^*(\lambda_{1})
< \alpha^*(\lambda_{2})$ follows.
\hfill $\Box$

\vspace{1pc}

\begin{lemma}
$\alpha^*(\lambda)$ 
is continuous 
for $\lambda \in ( \overline{\lambda}, \lambda^*)$.
\label{lem:conti-thresh-ld-alpha}
\end{lemma}
{\it Proof.}
Let $\lambda_{0} \in ( \overline{\lambda}, \lambda^*)$ be fixed.
Then
we have
$( \lambda_{0}, \alpha^*(\lambda_{0}) - \varepsilon ) \in \Gamma_{3}$ and
$( \lambda_{0}, \alpha^*(\lambda_{0}) + \varepsilon ) \in \Gamma_{1}$ 
for any chosen $\varepsilon \in (0, \alpha^*(\lambda_{0}))$.
By Lemma \ref{lem:openness-parameter-region},
there exists $\delta>0$ such that
$( \lambda, \alpha^*(\lambda_{0}) - \varepsilon ) \in \Gamma_{3}$ and
$( \lambda, \alpha^*(\lambda_{0}) + \varepsilon ) \in \Gamma_{1}$ 
hold for $\left\vert \lambda - \lambda_{0} \right\vert < \delta$.
Lemma \ref{lem:ex-thresh-ld-alpha} implies that 
$\alpha^*(\lambda_{0}) - \varepsilon 
< \alpha^*(\lambda) 
< \alpha^*(\lambda_{0}) + \varepsilon$,
and in turn 
$\left\vert \alpha^*(\lambda) - \alpha^*(\lambda_{0}) \right\vert
< \varepsilon$.
\hfill $\Box$

\vspace{1pc}

\begin{lemma}
We have $\lim_{\lambda \nearrow \lambda^*}\alpha^*(\lambda) = +\infty$.
\label{lem:asy-thresh-ld-alpha}
\end{lemma}
{\it Proof.}
Assume that
for $\lambda = \lambda^*$
there exists $\alpha_{0} > 0$ such that
$\gamma ( t; \lambda^*, \alpha_{0} ) \to ( p, 0 )$ as $t \to +\infty$,
where $\phi_{1} (\lambda^*) = \phi_{2} (\lambda^*) = p$.
We will show that this assumption leads us to a contradiction.
Indeed, in this case we have $\left( \lambda, \alpha \right) \in \Gamma_{1}$ 
for all $0 < \lambda < \lambda^*$ 
and 
$\alpha \geq \alpha_{0}$
by $\alpha^*(\lambda^*) \leq \alpha_{0}$ 
and 
Lemma \ref{lem:monotonicity-thresh-ld-alpha}.
Hence
$\gamma ( t; \lambda^*, \alpha ) \to ( p, 0 )$ as $t \to +\infty$ 
for all $\alpha \geq \alpha_{0}$.
Now,
again
through the transformation
\begin{equation*}
\left\{\begin{array}{ll}
t=-s, \\
U(s) = p - u(t), \\
V(s) = v(t),
\end{array}\right.
\label{eqn:transform*}
\end{equation*}
\eqref{eqn:dyn-sys} is equivalent to
\begin{equation*}
\frac{d}{ds}
\left[ \begin{array}{cc}
U\\
V
\end{array} \right]
= \left[ \begin{array}{cc}
V\\
\alpha V + H_{\lambda^*} \left( p - U \right) 
\end{array} \right]
\label{eqn:transformed-ODE*-hyperbolic}
\end{equation*}
for $s<0$ 
with $U(-\infty) = V(-\infty) = 0$.
We have one-parameter family of unstable manifold of $(U,V)=(0,0)$ with a branch
which enters
\begin{equation*}
W_{1}:= \left\{ \left( U,V \right) \in {\mathbb R}^2
\mid
U > 0, \ V > \Psi_{\alpha, \lambda^*}(U) \right\},
\label{eqn:convergence-region1}
\end{equation*}
\begin{equation*}
W_{2}:= \left\{ \left( U, V \right) \in {\mathbb R}^2
\mid
U > 0, \ 0 < V < \Psi_{\alpha, \lambda^*}(U) \right\}
\label{eqn:convergence-region2}
\end{equation*}
and 
approaches $(U, V) = (p, 0)$,
where
\begin{equation*}
\Psi_{\alpha, \lambda}(U)
\equiv -\frac{1}{\alpha}
H_{\lambda} \left( p - U \right).
\label{eqn:psi}
\end{equation*}
It is clear that
$V_{s}>0$ in $W_{1}$ and that $V_{s}<0$  in $W_{2}$.
Each of these branches of the unstable manifolds 
for $\alpha \geq \alpha_{0}$ 
is the graph of a function $V = \Phi (U; \lambda^*, \alpha)$
defined for $U \in [0, p]$.
Moreover,
similarly to \eqref{eqn:derivative-unstable-mfd},
we have
\begin{equation*}
\frac{d \Phi}{dU}(U; \lambda^*, \alpha)
= \alpha + \frac{H_{\lambda^*} \left( p - U \right)}{\Phi (U; \lambda^*, \alpha)}
\label{eqn:derivative-unstable-mfd-star1}
\end{equation*}
and 
\begin{equation*}
\frac{d \Phi}{dU}(0; \lambda^*, \alpha)
= \eta^+(\alpha)
= \alpha\ .
\label{eqn:derivative-unstable-mfd-star2}
\end{equation*}

\noindent As in Proposition \ref{prop:monotonicity-unstable-mfd-lambda},
the graph $V = \Phi (U; \lambda^*, \alpha_{2})$ stays 
above the graph $V =  \Phi (U; \lambda^*, \alpha_{1})$
for $\alpha_{0} \leq \alpha_{1} < \alpha_{2}$.
Let $\left( U_{0}, V_{0} \right)$ be the intersection point 
of the branch $V = \Phi (U; \lambda^*, \alpha_{0})$ 
of the unstable manifold
with the graph $V = \Psi_{\alpha_{0}, \lambda^*}(U)$.
Then
we have $U_{0}>0$,
$V_{0}>0$
and
\begin{equation*}
V_{0} = \Phi (U_{0}; \lambda^*, \alpha_{0}).
\label{eqn:intersection-point1}
\end{equation*}
By monotonicity and continuity,
there exists $\delta>0$ such that
\begin{equation*}
\frac{1}{2} V_{0} = \Phi (U_{0} - \delta; \lambda^*, \alpha_{0}).
\label{eqn:intersection-point2}
\end{equation*}
Let consider the value 
\begin{equation*}
\alpha > \max \left( \alpha_{0}, \frac{4 \lambda^*}{V_{0}} \right).
\label{eqn:parameter1}
\end{equation*}
Noting that $U_{0}<p$
and that
$0 > H_{\lambda^*} \left( p - U \right) 
\geq H_{\lambda^*} \left( 0 \right) 
=  - \lambda^*$
for $U \in [U_{0} - \delta, U_{0}]$,
we have
\begin{eqnarray*}
\frac{d \Phi}{dU}(U; \lambda^*, \alpha)
&=& \alpha + \frac{H_{\lambda^*} \left( p - U \right)}{\Phi (U; \lambda^*, \alpha)} \\
&\geq& \alpha + \frac{2 H_{\lambda^*} \left( p - U \right)}{V_{0}} \\
&>& \alpha - \frac{2 \lambda^*}{V_{0}} \\
&>& \frac{1}{2}\alpha
\label{eqn:derivative-unstable-mfd-star3}
\end{eqnarray*}
for $U \in [U_{0} - \delta, U_{0}]$. Apply the mean value theorem to obtain
\begin{eqnarray*}
\Phi (U; \lambda^*, \alpha)
&=& \Phi (U; \lambda^*, \alpha) 
- \Phi (U_{0} - \delta ; \lambda^*, \alpha) 
+ \Phi (U_{0} - \delta ; \lambda^*, \alpha) \\
&>& \frac{1}{2}\alpha \left\{ U - \left( U_{0} - \delta \right) \right\} 
+ \frac{1}{2} V_{0}
\label{eqn:unstable-mfd-star1}
\end{eqnarray*}
for $U \in [U_{0} - \delta, U_{0}]$
and in particular 
\begin{equation*}
\Phi (U_{0}; \lambda^*, \alpha)
> \frac{1}{2}\alpha \delta + \frac{1}{2} V_{0}.
\label{eqn:unstable-mfd-star2}
\end{equation*}
If necessary,
we can take $\alpha$ larger and satisfying
\begin{equation*}
\frac{1}{2}\alpha \delta + \frac{1}{2} V_{0}
> \frac{\lambda^*}{\alpha}
= \Psi_{\alpha, \lambda^*}(p)
\label{eqn:parameter2}
\end{equation*}
so that
\begin{equation*}
\Phi (U_{0}; \lambda^*, \alpha)
> \Psi_{\alpha, \lambda^*}(p).
\label{eqn:unstable-mfd-star3}
\end{equation*}
Finally we obtain
\begin{equation*}
\left\{ \left( U, V \right) \in {\mathbb R}^2
\mid V = \Phi (U; \lambda^*, \alpha), \ U_{0} < U < p \right\}
\subset W_{1},
\label{eqn:unstable-mfd-star4}
\end{equation*}
and thus 
\begin{equation*}
V_{0}
< \Phi (p; \lambda^*, \alpha)
= 0,
\label{eqn:unstable-mfd-star5}
\end{equation*}
which contradicts the fact $V_{0}>0$.
\hfill $\Box$

\vspace{1pc}

\noindent {\it Proof of Theorem \ref{thm:dynamical-threshold}.}
$\alpha^*(\lambda)$ is the desired threshold with respect to $\alpha$ 
by Lemma \ref{lem:ex-thresh-ld-alpha}.
Since $\alpha^*(\lambda)$ is strictly increasing and continuous
by Lemmas \ref{lem:monotonicity-thresh-ld-alpha} 
and 
\ref{lem:conti-thresh-ld-alpha} respectively,
the inverse function of $\alpha^*(\lambda)$, 
denoted by $\lambda \left( 0,0 \right) (\alpha)$, 
is well-defined
for $\alpha \geq 0$ and inherits the properties of monotonicity and continuity: Remark \ref{rem:ex-global-sol-origin}
and 
Lemma \ref{lem:asy-thresh-ld-alpha} complete the proof.
\hfill $\Box$

\section*{Acknowledgements}
This work was motivated and started
while the second named author was visiting the Riemann International School of Mathematics at the University of Insubria in Varese - Italy, during September 2017.
He would like to express his deepest gratitude for the warm hospitality.

\newpage 

\begin{flushleft}
{\it  Daniele Cassani\\
Dipartimento di Scienza e Alta Tecnologia\\
Universit\`{a} degli Studi dell'Insubria\\
and\\
RISM--Riemann International School of Mathematics\\
Villa Toeplitz, Via G.B. Vico, 46 -- 21100 Varese, Italy.\\
e-mail:} daniele.cassani@uninsubria.it
\end{flushleft}
\begin{flushleft}
{\it Tosiya Miyasita\\
Division of Mathematical Science,\\
Department of Science and Engineering,\\
Faculty of Science and Engineering,\\
Yamato University,\\
2-5-1, Katayama-cho, Suita-shi, Osaka, 564--0082, Japan\\
e-mail:} miyasita.t@yamato-u.ac.jp
\end{flushleft}


\begin{thebibliography}{99}

\bibitem{CFT}
D.~Cassani, L.~Fattorusso and A.~Tarsia, 
\textit{Nonlocal dynamic problems with singular nonlinearities and applications to MEMS}, 
Progress in Nonlinear Differential Equations and their Applications, 
{\bf 85},
Birkh\"{a}user (2014), 187–206.

 \bibitem{CKL}
 D.~Cassani, B.~Kaltenbacher and A.~Lorenzi, 
 \textit{Direct and inverse problems related to MEMS}, 
 Inverse Problems, 
 \textbf{25} 
 (2009), 
 105002 (22 pp).
 
\bibitem{MR1000727}
M. Chipot and F. B. Weissler,
{\it Some blowup results for a nonlinear parabolic equation 
with a gradient term},
SIAM J. Math. Anal.,
{\bf 20}
(1989),
886-907.

\bibitem{MR511133}
C. Conley,
{\it Isolated invariant sets and the Morse index},
CBMS Regional Conference Series in Mathematics,
{\bf 38},
American Mathematical Society, Providence, R.I.,
(1978).


\bibitem{MR3592652}
G. Flores,
{\it On the dynamic pull-in instability 
in a mass-spring model of electrostatically actuated MEMS devices},
J. Differential Equations,
{\bf 262}
(2017),
3597-3609.
 
 \bibitem{GK}
 J.-S.~Guo and N.I.~Kavallaris, \textit{On a nonlocal parabolic problem arising in electrostatic MEMS control}, Discrete Contin. Dyn. Syst. \textbf{32} (2012), 1723--1746.
 
\bibitem{MR2763360}
A. Haraux,
{\it Sharp decay estimates of the solutions to a class of nonlinear
second order ODE's},
Anal. Appl. (Singap.),
{\bf 9}
(2011),
49-69.

\bibitem{MR3380967}
A. Haraux and M. A. Jendoubi,
{\it The convergence problem for dissipative autonomous systems},
SpringerBriefs in Mathematics,
Springer, Cham; BCAM Basque Center for Applied Mathematics, Bilbao,
(2015).

\bibitem{KLNT}
 N.I.~Kavallaris, A.A.~Lacey, C.V.~Nikolopoulos and D.E.~Tzanetis, \textit{On the quenching behaviour of a semilinear wave equation modelling MEMS technology}, Discrete Contin. Dyn. Syst. \textbf{35} (2015), 1009--1037.


\bibitem{MR1270119}
L. A. Peletier, J. Serrin and H. Zou,
{\it Ground states of a quasilinear equation},
Differential Integral Equations,
{\bf 7}
(1994),
1063-1082.

\bibitem{MR1608005}
P. Souplet,
{\it Critical exponents, special large-time behavior and
oscillatory blow-up in nonlinear ODE's},
Differential Integral Equations,
{\bf 11}
(1998),
147-167.

\bibitem{VBM}
M.~Versaci, P.~Di Barba and F.C.~Morabito, \textit{Curvature-dependent electrostatic field as a principle for modelling membrane-based MEMS devices. A review}, Membranes \textbf{10} (2020), 1-51.   

\end{thebibliography}
\end{document}